\documentclass[12pt]{amsart}
\usepackage[paper=letterpaper, margin = 3cm]{geometry}
\usepackage{graphicx} 
\usepackage{amsfonts}
\usepackage{mathtools}
\usepackage{amsmath}
\usepackage{amssymb}
\usepackage{amsthm}
\usepackage{tikz}
\usetikzlibrary{positioning, arrows}
\usepackage{parskip}

\usepackage{comment}

\newtheorem{theorem}{Theorem}[section]
\newtheorem{lemma}[theorem]{Lemma}
\newtheorem{proposition}[theorem]{Proposition}
\newtheorem{corollary}{Corollary}[theorem]
\newtheorem{remark}{Remark}[theorem]

\theoremstyle{definition}
\newtheorem{definition}[theorem]{Definition}
\newtheorem{example}{Example}[theorem]

\numberwithin{equation}{section}

\newcommand{\m}[1]{\bigl(\begin{smallmatrix}#1\end{smallmatrix}\bigr)}
\newcommand{\M}[1]{\begin{pmatrix}#1\end{pmatrix}}

\newcommand{\mycomment}[1]{}

\title{Freedom is all you need}
\author{Julian Bushelli, James Eldred Pascoe}
\date{November 2026}

\begin{document}

\begin{abstract}
    
The functional calculus allows for the rigorous analysis of functions of matrices or even operators. The free functional calculus gives the natural class of functions over any category. We initiate the representation theory for taking a generic free functional calculus to the classical nc functional calculus going back to Taylor, Voiculescu and many others, with an eye towards learning of categories. The functional calculus on direct sum closed and locally similarity (i.e. change of coordinates) invariant domains inherits the niceties of complex analysis and large-scale parallel networks, such as automatic analyticity, free universal monodromy, which gives a hypothetical explanation for beneficial overfitting. 
\end{abstract}
\maketitle

\section{Introduction}

It has been long understood that if a machine learning algorithm is allowed to optimize over too many parameters, it will ``over-fit" to the training data and fail to properly model the underlying phenomena. However, it has recently been noticed that this does not happen to computational neural networks, but instead increasing parameters continues to improve the predictive power of the model indefinitely [ZBHRO21]. 

We suggest a theoretical explanation of this phenomenon; namely, that a computational neural network has the structure of a ``free function," and thus benefits from the many nice properties of these functions. Crucially, free functions respect direct sum structures, which makes them particularly suited to understanding parallel systems. These free functions can be considered as a non-commutative analogue of complex analytic functions, and satisfy even stronger theorems (due to their dimension-free nature) than typical complex analytic functions.  This theory is well studied when the functions take matrix or operator variables [KVV14, JMS21, HKM12, Tay72], which in turn has many roots in control theory and the study of system-flow diagrams; the central theory extends elegantly to the case of so-called ``additive categories," [B25] but while there is reason to believe that the theory is even more broad-reaching, these limits of the theory are not well understood [Pas23]. 

Neural networks fail to satisfy the requirements of an additive category. So as to easily apply these results, we will ``additivize" this category by representing it as a category of Hilbert spaces. The topological spaces will together be given a cross-compatible set of reproducing kernel Hilbert spaces, one on each space, so that the neural networks can be represented as linear composition operators between these vector spaces. This can be thought of as a ``computable model" of the underlying topological space. The theory of RKHSs has been leveraged for many purposes, including statistical analysis, machine learning, control theory, signal processing, and function theory, as well as others. 

A main question of this paper is the following; if we can use RKHSs to create is an observable model of a family of topological spaces, when will this model faithfully represent the underlying phenomenon, and when will the complexity of the phenomenon preclude this possibility? We give a partial answer; any ``sufficiently decomposable" family of second-countable Tychonoff spaces can be faithfully represented. If the spaces are not second-countable, the faithfulness of the observable model is no longer guaranteed. This is all developed within section III.

Section IV provides generalizations of several results in free analysis that are useful for analyzing neural nets, most notably an extension of the notion of analytic continuation. This theory is predicated upon an additive structure, and provides a strong analogue of complex analysis for such spaces. 

Section V then discusses the specific instance of neural nets, giving an explanation for ``beneficial overfitting" by leveraging the constructed additive structure together with the free functional calculus.  

\section{Neural nets, topologies, and categories}

\subsection{Topologies and categories}
Suppose we have a collection of topological spaces and continuous maps. Given continuous maps $f: X \rightarrow Y$ and $g: Y \rightarrow Z$ we can compose them to define the continuous map $g \circ f: X\rightarrow Z$. 
\begin{center} \begin{tikzpicture}[node distance=1cm, auto]
  \node (C) {$Y$};
  \node(b) [below of = C]{};
  \node (Q) [left of=b, xshift = -1cm] {$Z$};
  \node (R) [right of=b, xshift = 1cm] {$X$};
  \draw[->] (C) to node [swap] {$g$} (Q);
  \draw[->] (R) to [name = fx, swap] node {$f$} (C);
  \draw[->] (R) to node [swap] {$g\circ f$} (Q);
\end{tikzpicture}\end{center}
If we allow all such compositions as well as the trivial identity functions, then this collection of spaces and maps is called a \textit{category of topological spaces.} This is an archetypal example of a category, and is a worthwhile first example in many circumstances. In fact, every small category can be faithfully represented as a category of (discrete) topological spaces, due to the Yoneda lemma [Rie17]. Thus, one may use trivial kernels to model any category.

\begin{definition}[Category]
    A \textbf{category} $C$ consists of the following 
    \begin{itemize}
    \item A collection of \textbf{objects}, denoted $ob(C)$
    \item For every $u, v\in ob(C)$, a collection of \textbf{morphisms} (also called \textit{arrows} or \textit{maps}) from $u$ to $v$, denoted $mor(u, v)$. 
    \begin{itemize}
        \item If $x \in mor(u, v)$ and $y\in mor(v, w)$ then there is a morphism $yx := y\circ x \in mor(u, w)$; this composition must be associative. 
        \item There exists an element $1 \in mor(u) := mor(u,u)$ called the \textbf{identity} so that $1x = x$ and $x1 =x$ for any $x$ with appropriate domain and codomain. 
    \end{itemize}
    \end{itemize} 

    The notation $x: u \rightarrow v$ is equivalent to $x \in mor(u,v)$. 
\end{definition}

As discussed above, the most important examples are categories whose objects are topological spaces and whose morphisms are continuous maps. 

We must also define representations of categories. 

\begin{definition}
    Let $C, D$ be categories. A \textbf{functor} $f: C\rightarrow D$ is a mapping between categories satisfying the following. 
    
    For all $u, v, w \in ob(C)$, and for all $x\in mor(u, v)$ and $y\in mor(v, w)$:  
    \begin{itemize}
        \item There is an object $f(u) \in ob(D)$.
        \item There is a morphism $f(x) \in mor(f(u), f(v))$ for each $x \in mor(u, v)$.
        \item For the identity $1_{u} \in mor(u)$, the image $f(1) = 1_{f(u)} \in mor(f(u))$ is the identity. 
        \item The functor respects composition, that is $f(yx) = f(y)f(x)$.
    \end{itemize}

    The above functor is sometimes said to be ``covariant." If the last relation above is replaced with $f(xy) = f(y)f(x)$, the mapping is called a  \textit{contravariant functor}. 

    A functor is said to be \textbf{faithful} if for every pair of objects $u, v$ and all morphisms $x,y: u\rightarrow v$ we have that \[f(x) \cong f(y) \iff x \cong y.\] 
\end{definition}

\begin{example}[Dualization]
    An important example of a contravariant functor is the ``dual map" which sends a vector space $X$ to the space of linear functionals $X^*$, and sends a linear map $T: X \rightarrow Y$ to its dual map $T^*: Y^*\rightarrow X^*.$
    
    The inclusion of a vector space $X \subset X^{**}$ in its double dual is a good example of a faithful functor. 
\end{example}
\begin{example}[Forgetful functors]
Another important example of a faithful functor is the so-called ``forgetful" functor that takes vector spaces and ``forgets" their vector space structure, sending them to their underlying topological spaces and sending the linear transformations to continuous maps between those topological spaces. 
\end{example}
\begin{example}[Composition operators]
One can map a topological space $X$ to the space $C(X)$ of continuous (not necessarily linear) functionals on $X.$ Then one can send a continuous map $\varphi: X\rightarrow Y$ to the composition operator (i.e. pullbacks) $C_{\varphi}: C(Y) \rightarrow C(X).$ This is a contravariant functor, and functors of this sort will be discussed at length later. 
\end{example}

It is possible to get some idea of how compatible two functors $F, G: C \rightarrow D$ are by looking at the relations within the target category. 

\begin{definition}
    Consider two functors $F, G: C\rightarrow D$. A collection of morphisms $\eta$ consisting of morphisms $\eta_u: F(u) \rightarrow G(u)$ within $D$ for each $u\in obC$ is called a \textit{natural transformation} from $F$ to $G$ if the following square commutes for every $u, v \in ob C$, $a :u\rightarrow v.$ 
     \begin{center} \begin{tikzpicture}[node distance=2cm, auto]
  \node (X) {$F(u)$};
  \node (fx) [right of=X] {$F(v)$};
   \node (Y) [below of=X] {$G(u)$};
  \node (fy) [right of=Y] {$G(v)$};
  \draw[->] (X) to node [swap] {$\eta_u$} (Y);
  \draw[->] (X) to node {$F(a)$} (fx);
  \draw[->] (Y) to node {$G(a)$} (fy);
  \draw[->] (fx) to node [swap] {$\eta_v$} (fy);
\end{tikzpicture}\end{center} which is interpreted algebraically as \[\eta_{cod(a)}F(a) = G(a)\eta_{dom(a)}\]
for each morphism $a$ of $C.$
\end{definition}

We then define the category $D^C$, called the \textit{functor category}, whose objects are functors $F: C\rightarrow D$ and morphisms are natural transformations $\eta: F_1 \rightarrow F_2$. 

Classically the notion of a natural transformation shows up within algebra as inner automorphisms. In analysis we have intertwining as a simple example of natural transformations, with the solutions to the Sylvester equation $AX - XB =0$ telling us when certain natural transformations exist. 

\subsection{Neural Networks}

A (computational) neural net should be considered a special case of a continuous map. Let us consider a neural net taking $n$ real inputs and giving us $m$ real outputs. 

\begin{definition}[Neural Networks]
    Define the function $ReLu: \mathbb R\rightarrow \mathbb R^+$ by $ReLu(x) = \max\{0,x\}$.
    
    A depth-zero (computational) neural network $F: \mathbb R^n \rightarrow \mathbb R^m$ is a function of the form \begin{equation}F(x) = (F_1(x), \dots, F_m(x)) = (ReLu\left(\sum_{i \in [n]}a_{i1}x_i\right), \dots, ReLu\left(\sum_{i\in [n]}a_{im}x_i\right)).\end{equation}

    A depth $\ell$ neural network is a composition of $\ell$ depth-zero neural nets. This number corresponds to the number of ``hidden layers" in the neural net. 
\end{definition}

\begin{remark}
    There are several common variants from what is described above. Often  $ReLu$ is replaced by a smooth approximant, such as $exp$ or a sigmoid function. Similarly, often the last layer of the neural network will not have $ReLu$ applied, so that the function is not always positive. 

    Since the following discussion does not rely on the details, we will only consider the ``standard" case, defined above. 
\end{remark}

The category of neural nets $CNN$ is naturally defined in the following way. 

\begin{itemize}
    \item The objects of $CNN$ are vector spaces $\mathbb R^n$, $n\in \mathbb N$.
    \item The morphisms are neural nets of any depth; 
    \item The identity morphism on $\mathbb R^{n}$ is the neural net \[F(x) = (x_1, \dots, x_n);\]
    \item Composition of morphisms is function composition of neural nets; 
    \begin{itemize}
        \item If $f$ has depth $n$ and $g$ has depth $m$, $f\circ g$ has depth $n+m+1$.
    \end{itemize}
\end{itemize}

Note that, as constructed, the range of these functions is only the positive orthant of each $\mathbb R^n$. Thus, it would be very reasonable to view the objects as only being the positive orthant $\mathbb R^{n+}$. 

There is an addition operation for $ f, g \in mor(\mathbb R^n, \mathbb R^m)$, defined by $(f+g)(x) := f(x)+g(x) \in \mathbb R^m$. We might hope that this would allow us to form an additive category and apply the free functional calculus to these maps; indeed, it is clear from the definition of addition that $(f+g)\circ h = f\circ h+g\circ h,$ so left-distributivity holds. However, since our functions are not linear, in general we see that \[h\circ (f+g) \neq h\circ f+h\circ g.\]

For example, if $h: \mathbb R^2\rightarrow \mathbb R$ is defined as $h(x, y) = ReLU(x-2y)$ and $f(x_0) = (0,3), g(x_0) = (1, 0)$ then \[\begin{split}h\circ (f+g)(x_0) &= h(1,3)= ReLU(1-6) = 0\\ (h\circ f + h\circ g)(x_0) &= h(1,0)+h(0,3) = ReLU(1)+ ReLU(-6)  = 1
\end{split}\]

Thus the category of neural nets is therefore not an additive category, though it is a \textit{Cartesian left-additive category}, which is a type of category that has been studied in other contexts. In these categories matrix computations do not work normally, which severely limits the analytic techniques at our disposal. 

In order to regain these techniques, we must find a way to represent this category of topological spaces as an additive category, by way of composition operators on reproducing kernel Hilbert spaces. 

\subsection{Reproducing Kernel Hilbert spaces}

Here we will give a basic summary of reproducing kernel Hilbert spaces (RKHS), and show how it is a good example of a contravariant functor. 

A \textbf{reproducing kernel Hilbert space} $\mathcal H_k$ is a Hilbert space of functions $f: \Omega \rightarrow \mathbb R$ (or $\mathbb C$)  so that point evaluation defines a linear functional, that is $\lambda_xf: \mathcal H_k \rightarrow \mathbb R$ defined by $\lambda_x f = f(x)$ is a linear functional. 

It is called a reproducing kernel Hilbert space because there is a bijection between such spaces and positive semidefinite binary ``kernel" functions. 
\begin{definition}
    A function $\mathcal K: \Omega \times \Omega \rightarrow \mathbb R$ is called a \textit{reproducing kernel} if it satisfies $\mathcal K(x, y) = \mathcal K(y, x)$ and the inequalities \[\sum_{i,j}^n c_ic_j \mathcal K(x_i, y_j) \geq 0.\] 
\end{definition}

Due to the Moore-Aronszajn theorem, we are able to construct a unique RKHS associated with each reproducing kernel. The functions in this space are generated by taking the closure of the linear combinations of point functions $K_x := K(x, -)$ and the inner product is given by defining $\langle K_x, K_y\rangle = K(x, y)$ and expanding any function into its components. So if $f = \sum_i c_i\mathcal K_i(x_i, -),$ \[\langle f, f\rangle = \sum_{i,j} c_ic_j \mathcal K(x_i, y_j) = \|f\|^2.\] This inner product is well-defined if $\mathcal K$ is \textit{strictly} positive definite (that is, the inequality above is strict), but we must mod out the null elements otherwise. This shows us that the point functions $k_x$ are all linearly independent if and only if the kernel is strictly positive definite. 

Note that in the sequel we will follow the conventions of RKHS literature in saying ``strictly positive definite" when we mean the inequality is strict, and treating the terms ``positive definite" and ``positive semi-definite" as synonymous. 

\begin{example} The following are examples of reproducing kernels. 
\begin{itemize}
    \item the linear kernel $k(x,y) = y^Tx$. 
    \item the degree-$d$ polynomial kernel $k(x, y) = (a y^Tx+c)^d$, with $a, c > 0$. These give the degree $d$ polynomials. 
    \item the Drury-Arveson kernel, $k(x, y) = \frac{1}{1-\langle x, y\rangle}$, which generates the Drury-Arveson space, a multivariable analogue of the the Hardy space $H^2$. 
    \item The Gaussian kernel (RBF) on a metric space, where $k(x, y) = \exp(-\gamma^2 d(x, y))$
\end{itemize}
\end{example}

\begin{definition}Two kernels $k, h$ are said to be \textbf{unitarily equivalent} if $h(x,y) = w(x)w(y)k(x,y)$ for some continuous nonzero function $w$. \end{definition}

\section{Additivizing Topological Categories}

In an approach popularized by Koopman [Koo31] and used by many others, we will linearize our space to create an additive category. The following line of reasoning works in any category of topological spaces with continuous functions as morphisms. We shall deal in generality, then specialize to the category of neural nets for more concrete results.

\subsection{Composition operators}
\begin{definition}
    On each topological space $X$, we can construct a reproducing kernel Hilbert space of continuous functionals $\mathcal H_X$, with elements $f: X \rightarrow \mathbb R$. \footnote{this is essentially slicing the category of vector spaces over the reals}

    Each function $\varphi: Y \rightarrow X$ defines a linear operator; we define $C_\varphi: \mathcal H_X \rightarrow C(Y)$ by $C_\varphi f(y)= f(\varphi (y))$ for each $y \in Y$, where $C(Y)$ is the set of continuous functionals on $Y$. 
\end{definition}

It is straightforward to check that $C_\varphi$ is linear. Indeed, \[C_\varphi(f+\alpha g)(y) = (f+\alpha g)(\varphi (y)) = f(\varphi(y))+\alpha g(\varphi(y))= C_\varphi f(y)+\alpha C_\varphi g(y).\] There is a simple addition structure, with $C_\varphi + K_\psi$ defined by \[(C_\varphi + C_\psi)f = C_\varphi f+C_\psi f.\] Putting these two facts together, it is simple to check that this addition distributes on both sides. Defining the additive inverse as $-C_\varphi$ sending $f \mapsto -f\circ \varphi$, we see that there is a vector space of composition operators generated by the maps $\varphi: Y\rightarrow X$. 

The goal is to give RKHSs to each space to form a pre-additive category, and then include $\mathcal H_X \oplus \mathcal H_Y$ as well as all other such direct sums to create an additive category. Theorem 2.19 of [SaSa16] tells us that these direct sums are themselves isomorphic to an RKHSs, so our additive category consists of RKHSs. 

There is a problem, however; if we have an RKHS $\mathcal H_Y$ of functions on the space $Y$, there is no guarantee that $C_\varphi f \in \mathcal H_Y$ if $f \in \mathcal H_X$. In particular, it is easy to find examples of kernels on $X, Y, X \times Y$ so that $\mathcal H(X) \not \subset H(X\times Y)$, or more to the point $H(X) \oplus H(Y) \not \subset H(X\times Y)$. For example, if $k_X$ is the degree 2 polynomial kernel, and $k_{X\times X}$ the degree $4$ polynomial kernel, the restricted space will not match. 

In order to deal with these subtleties, we introduce the notion of a ``universal kernel family," which is a family of kernels on a category $T$ of topological spaces which respects products. 

\subsection{Universal families of kernels and faithful representations}

\begin{definition}
    Let $T$ be a cartesian category of topological spaces, which is to say it contains all finite products. For each $T_i \in obT$ choose a point $0_i =: 0$. Let $k$ be a family of kernels, with a unique kernel $k_X$ for each topological space $X \in obT$. 

    Define \[\mathcal H_{k}(X\times Y)|_{(X, 0)}\] as the RKHS generated by restricting $k_{X\times Y}: (X\times Y)\times (X\times Y) \rightarrow \mathbb R$ onto $(X, 0)\times (X, 0).$

    We say that $k$ is a $\mathcal T -$universal family of kernels if it satisfies \[\mathcal H_{k}(X\times Y)|_{(X, 0)} \cong \mathcal H_k(X),\] where the isomorphism is an unitary equivalence of RKHSs, as well as \[\mathcal H_{k}(X\times Y)|_{(0, Y)} \cong \mathcal H_k(Y).\]
    
    This means \[k_{X \times Y}((x_1, 0), (y_1,0)) = w(x_1)w(x_2)k_X(x_1, x_2)\] for some continuous nonzero function $w$. 

    Where $\mathcal H_{k}(X) := \mathcal H_{k_X}(X)$, as these kernels are unique. 

    We will also sometimes use the notation $\mathcal H_X := \mathcal H_{k_{X}}(X)$, as it can only mean one thing for given $T, k$. 
\end{definition}

\begin{example}
    If $T$ is a category of metric spaces with continuous functions, with metric $d_i$ on space $T_i$, then the Gaussian kernel \[e^{-\gamma d(x,y)}\] for $\gamma >0$ defines a $T$-universal family of kernels.
\end{example}

\begin{example}
    The polynomial kernels \[k(x, y) = (ax^Ty + c)^d,\] where $x, y \in \Omega \subset \mathbb R^n$, where $a, c, d$ are fixed and $n$ depends on the topological space in question. These define a universal family of kernels over the category of finite dimensional real vector spaces, or categories of subsets thereof. The finite-dimensional space $\mathcal H_{\mathbb R^n}$ is the space of $n-$variate real polynomials of degree $d$, as shown in [Zwi08] and [EEM22].

\end{example}

\begin{proposition}
    If we choose a universal family of kernels $k$, then for topological spaces $X, Y \in obT$ the induced RKH spaces satisfy $H_X \oplus H_Y \subset H_{X\times Y}$.
    
    In particular, $H_X \oplus H_Y$ is itself an RKHS when $\oplus$ is considered as the internal direct sum given by addition. 
\end{proposition}

\begin{proof}
    Consider the inclusions  $X \hookrightarrow X \coprod Y \hookrightarrow X \times Y$. This can be realized by identifying $X$ with $(X, 0)$, and $Y$ with $(0, Y)$. 

    Then by [Theorem 2.19, 3], we can create a RKHS isomorphic to $\mathcal H_X \oplus \mathcal H_Y$ on $X \coprod Y$ by simple addition of the functions. 

    It is clear from the definition of a universal family of kernels that $\mathcal H_X \hookrightarrow \mathcal H_{X \times Y}$, and similarly for $Y$. We realize this inclusion as $\mathcal H_X \cong \mathcal H_{X \times Y}|_{(X, 0)}$, and similarly for $Y$. 

    We then identify $\mathcal H_{X} \oplus \mathcal H_Y$ with the internal direct sum $\mathcal H_{X \times Y}|_{(X, 0)} \oplus \mathcal H_{X \times Y}|_{(0, Y)} \subset \mathcal H_{X \times Y}.$
\end{proof}

We have discussed at length the properties of the reproducing kernel Hilbert spaces themselves, but we wish to apply a morphism functional calculus to the modified category. Are these composition operators truly linear maps in this category of RKHSs? Is this mapping faithful? We will address these questions.

\begin{proposition}\label{between spaces}
    If $K$ is a $T-$universal family of kernels, then $f: X \rightarrow Y$ becomes a composition operator $C_f: K_Y \rightarrow K_X$. 
\end{proposition}
\begin{proof}
    Because $K$ is a universal family of kernels, we shift our focus to $K_{X\times Y}$. In this context, $f$ is an endomorphism defined by mapping $(x, y) \mapsto (0, f(x))$; in particular $(x, 0) \mapsto (0, f(x))$. This map is continuous, as $f^{-1}(0, U) = (f^{-1}(U), Y)$. 

    We now denote the kernel function on $K_Y$ as $k'$, and on $K_X$ as $\tilde k$.

    Now let $k'_y \in K_Y$. We extend this to $k_{(0, y)} \in K_{X\times Y}$. Then \[ \begin{split}\langle \tilde k_x, C_fk'_y \rangle &=\langle k_{(x, 0)} ,C_fk_{(0, y)}\rangle\\& = k_{(0, y)}(f(x, 0))\\&=k_{(0, y)}(0, f(x))\\&= k'(y, f(x))\end{split}\]
    and so $C_fk'_y \in K_X$.
\end{proof}

This can perhaps be seen more easily by looking at the adjoint: $\langle C^*_f\tilde k_x, k'_y\rangle = \langle \tilde k_{f(x)}, k'_y \rangle =k'(y, f(x))$, so $C_f^*\tilde k_x \in K_Y$. 

This shows us that our additive category representation sends continuous functions to bona fide linear maps between our RKHSs.

Next is given to show that the mapping of morphisms is faithful. Recall that if a kernel is (strictly) positive definite, then each $k_x$ is distinct for distinct $x$, and so the mapping is injective. Indeed, if ever we have that $k(x, y) = k(y,y)=k(x, x)$, then if we set $c_x = 1$, $c_y = -1$ we have 
\[k(x,x) -k(y, x) -k(x, y) + k(y, y) = 0.\]
Thus we can satisfy the following criterion by choosing a positive definite kernel, rather than one that is positive semi-definite.

\begin{proposition}
    If the maps $x\mapsto k_x$ are injective for the spaces $X_1, X_2$, then each distinct continuous $\varphi: X_1 \rightarrow X_2$ yields a distinct composition operator $C_\varphi$. 
\end{proposition}
\begin{proof}
     We can consider $f: X\rightarrow Y$ as a continuous function $\tilde f: X\times Y \rightarrow X\times Y$ as in Proposition \ref{between spaces}, so without loss of generality we assume any function $f : X\rightarrow X$ for some space $X$.
    
    Consider $f$, a continuous function on $X$, and let $ C_f$ be its composition map. We define the adjoint $C^*_f$ as satisfying $\langle k_x, C_fk_y\rangle = \langle C^*_fk_x, k_y\rangle$.  Then \[\begin{split}
        \langle k_x, C_fk_y\rangle &= C_fk_y(x) \\ &= k_y(f(x))\\&=k(y, f(x)) \\&=k_{f(x)}(y).
    \end{split}\]

    On the other hand, $
        \langle C^*_fk_x, k_y\rangle = C^*_fk_x(y),$
        and therefore \[C^*_fk_x = k_{f(x)}.\] 
        
        Now let $f \neq g$, that is $f(x) \neq g(x)$ for some $x \in X$. Then \[C^*_fk_x \neq C^*_gk_x,\]
        in particular there exists $y\in X$ so that $k_{f(x)}(y) \neq k_{g(x)}(y).$

        This means that \[\langle k_x, C_fk_y\rangle = \langle C^*_fk_x, k_y\rangle \neq \langle C^*_gk_x, k_y\rangle = \langle k_x, C_gk_y\rangle.\]
\end{proof}

These aforesaid results can be largely summarized in the following way.

\begin{theorem}
    Given a cartesian category $\mathcal T$ of topological spaces, a $\mathcal T$-universal family of kernels $\mathfrak K$ induces a covariant functor $\mathcal K$ which maps $\mathcal T$ to an additive category of reproducing kernel Hilbert spaces, mapping each continuous function $\varphi$ to $C_\varphi^*$. 

    This functor can also be considered as a contravariant functor, mapping $\varphi \mapsto C_\varphi$. 

    Additionally, $\mathcal K(\mathcal T)$ has the following superadditive property: \[\mathcal K(X) \oplus \mathcal K(Y) \subset \mathcal K(X\times Y).\]

    The functor $\mathcal K$ is faithful if the family of kernels $\mathfrak K$ is strictly positive definite. 
\end{theorem}

\subsection{Universal kernels in second countable spaces}

We now turn our attention to when we can construct such kernels. It is already clear that on any category of metric spaces we can construct a universal family of kernels, simply by using the Gaussian kernel. 

\begin{definition}
    A category of topological spaces $\mathcal T$ is \textbf{decomposable} if every sequence of projections terminates. 

    That is, every object in the category is generated by a product of spaces which are not themselves products of any other spaces in $T$. These indecomposable spaces are herein called \textbf{atoms,} or \textbf{atomic spaces.}
\end{definition}

\begin{theorem}
\label{Tychonoff}
    Any decomposable category of second countable Tychonoff spaces can be equipped with a universal family of continuous and strictly positive definite kernels.  
\end{theorem}

To prove this, we will appeal to the following two classical results.

\begin{proposition}[Theorem 32.1 in Mun00]\label{normal}
    Every second countable regular topological space is normal. 
\end{proposition}

\begin{proposition}[Urysohn's Lemma, 33.1 in Mun00]
    A topological space $T$ is normal if and only if any two disjoint closed sets $A$ and $B$ there exists a function $f: T\rightarrow \mathbb R$ so that $f(A) = 1$ and $f(B) = 0$.
\end{proposition}

These two classical propositions are used primarily in the following lemma. 

\begin{lemma}
    If $T$ is a second countable normal space, then for every open set $U \subset T$ there is a function $f: T \rightarrow \mathbb R$ so that $1 >f(U) > 0$ and $f(T\setminus U) = 0$. 
\end{lemma}
\begin{proof}
    For each open set $U$ we do the following. Let $B = T\setminus U$ be a closed subset of $T$, and let $\mathcal B$ be a countable basis of the topology. Then because $U$ is an open set, for each $a \in U$ there is a $a \in U_a \subset U$ such that $\overline U_a \subset U$. Indeed there must be such an open set in $\mathcal B$. Thus we have a countable family of such $U_a$.

    Then since $T$ is normal, we can create functions $f_i$ so that $f_i|_{U_i} = 1$, $f_i|_{X \setminus B} = 0$. Then choosing some nonzero $\alpha_i$ so that $\sum_i^{\infty}\alpha_i = 1$, we define $f_U = \sum_{i}\alpha_if_i$. This is nonzero everywhere on $U$, and zero everywhere else. 
\end{proof}

We require the following several lemmata to complete the argument. 

\begin{lemma}
\label{regular kernel}
    Every second countable regular Hausdorff space has a strictly positive definite continuous kernel.

    In particular, there is such a kernel of the form $k(x,y) = \sum_i\alpha_if_i(x)f_i(y)$ where $\sum_i\alpha_i = 1$. 
\end{lemma}
\begin{proof}
    Let $T$ be a second countable regular Hausdorff space; by Proposition \ref{normal}, T is normal. We can choose a countable open covering $U$ of $T$ so that every open set of $T$ is a union of $U_i \in U$; call this countable basis $\mathcal B$.
    By the previous lemma, for each open set $U_i \in \mathcal B$ there is a function $f_i$ so that $1>f_i|_U > 0$, but $f_i|_{T\setminus U} = 0$. 

    Choose some weights $\alpha_i$ so that $\sum_i\alpha_i^2 = 1$. We then define a kernel on $T$ by \[k(x,y) = \sum_i\alpha_i^2f_i(x)f_i(y),\]
    which is continuous and always bounded by $1$. 

    Now it remains to show that $k(x, y)$ is strictly positive definite. We do this by showing that the evaluation functionals $k_x = k(x, -)$ are linearly independent. Thus we choose a finite set of distinct points $x_1, \dots, x_n \in T$. Because $T$ separates points, an intersection of open sets guarantees an open set containing $x_i$ but missing all the others, which in turn guarantees one such open set in $U$. 
    
    We consider the so-called ``feature map" representation of $k_x$, namely considering $k_x$ as a vector in $\mathbb R^\mathbb N$ with the $i$-th coordinate of $k_x$ given by $f_i(x)$, and $\langle k_x, k_y\rangle_{\mathbb R^\mathbb N} = k(x, y)$. But the previous paragraph shows that any finite set of points will give a linearly independent set of vectors under the feature map.

    Then \[\sum_i\sum_ja_ia_jk(x_i, x_j) = \sum_i\sum_ja_ia_j\langle k_{x_i}, k_{x_j}\rangle = \|\sum_ia_ik_{x_i}\|^2\]

    which is strictly positive if $x_i$ are distinct and $a_i$ not all zero. 
    
\end{proof}

\begin{lemma}\label{products}
    If $X, Y$ are second countable regular Hausdorff spaces with kernels as in the above lemma, then \[k(x,y) = \sum_{i,j}\alpha_i\beta_j(f_i\times g_j)(x)\cdot(f_i\times g_j)(y)\] is a strictly positive definite continuous kernel. 

    In particular, the restriction of this kernel to $(X, 0)$ is equal to $\gamma k_X$ for some scalar $\gamma > 0$.
\end{lemma}
\begin{proof}

    Consider the product topology on $X \times Y$, where $X$ has countable open generating set $U$, and $Y$ has $V$. Then if the open set $\bigcup_{i} U_{i} =u \subset X$ and $\bigcup_{k}V_{k} = v \subset Y$ for $\{U_i\} \subset U, \{V_i\}\subset V$ then \[\bigcup_{i,j}U_i \times V_j = u\times v \subset X\times Y.\]

    So $U \times V$ is a countable open generating set of the topology on $X\times Y$.

    Then for $f: X \rightarrow \mathbb R$, $g: Y \rightarrow \mathbb R$ we define $f \times g: X\times Y \rightarrow \mathbb R$ by $(f\times g)(x, y) = f(x)g(y) \in \mathbb R$ for $(x, y) \in X\times Y$. 

    Since $\sum_i\alpha_i=1=\sum_j\beta_j$, we also have that $\sum_i\sum_j\alpha_i\beta_j=1$. Then using the kernels $k_X(x_1, x_2) = \sum_i\alpha_if_i(x_1)f_i(x_2)$ and $k_Y(y_1, y_2) = \sum_i\beta_ig_i(y_1)g_i(y_2)$ we can define the kernel \[\begin{split}k_{X\times Y}(z_1, z_2) &= \sum_{i,j}\alpha_i\beta_j(f_i\times g_j)(z_1)\cdot(f_i\times g_j)(z_2) \\&= \sum_{i,j}\alpha_i\beta_jf_i(x_1)\cdot f_i(x_2)\cdot g_j(y_1)\cdot g_j(y_2)\end{split}\]
    for $z_1, z_2 \in X\times Y.$ But this has exactly the structure we constructed in the previous lemma, so this is a strictly positive definite continuous kernel. 

    All that remains now is to show that it restricts appropriately. Thus let us restrict to $X$ by evaluating at $(x, 0)$ for the chosen points $0$ of $Y$. Then for some $g_j$ we have $g_j(0) \neq 0$, and so \[\begin{split}k_{X\times Y}((x_1,0), (x_2, 0)) &= \sum_{i,j}\alpha_i\beta_jf_i(x_1)\cdot f_i(x_2)\cdot g_j(0)\cdot g_j(0)\\&=\sum_j\sum_i\alpha_if_i(x_1)f_i(x_2)\cdot \beta_jg_j(0)^2\\&=\sum_j\beta_jg_j(0)^2\sum_i\alpha_if_i(x_1)f_i(x_2)\\&= \gamma \sum_i\alpha_if_i(x_1)f_i(x_2) \leq 1.\end{split}\]
\end{proof}

This next lemma has the flavor of a ``converse" lemma to the above.

\begin{lemma}
\label{restrictions}
    If $T$ is a space with a strictly positive definite continuous kernel on it of the form $k(x,y) = \sum_i\alpha_if_i(x)f_i(y)$ where $\sum_i\alpha_i = 1$, then any restriction to a subspace has a strictly positive definite kernel of the same form. 

    In particular, if there are kernels with the appropriate restriction property on $X\times Y\times X\times Z\times W$, $X\times Y$, and $X$, then it can be extended to $X\times Z$ in a way that makes the following diagram commute. 

    \begin{center}
        \begin{tikzpicture}[node distance=1.5cm, auto]
    \node (xyxzw) {$X\times Y\times X\times Z\times W$};
\node (xy) [below of = xyxzw, left of = xyxzw] {$X\times Y$};
\node (xz) [node distance = 3cm, right of = xy] {$X\times Z$};
\node (x) [below of = xy, right of = xy] {$X$};

\draw[->] (xyxzw) to node {$\pi_1$} (xy);
\draw[->] (xyxzw) to node {$\pi_2$} (xz);
\draw[->] (xy) to node {$\pi_{11}$} (x);
\draw[->] (xz) to node {$\pi_{21}$} (x);

\end{tikzpicture}
    \end{center}
    In particular we require the relation \[k_{X\times Y\times X\times Z\times W}|_{(X, 0, 0,0, 0)} \cong k_{X\times Y\times X\times Z\times W}|_{(0, 0, X,0, 0)} \cong k_X,\] that is it behaves appropriately along both restrictions to $X$.

\end{lemma}
\begin{proof}
    Let $X \subset T$ be a subspace. Then via restriction to the subspace topology, we obtain a countable basis $\mathcal B_X$ from the countable basis $\mathcal B$ of $T$. 

    Similarly, for each open set $b \in \mathcal B_X$ we have a continuous function $f_b$ that is strictly positive on $b$ and zero elsewhere, obtained by restricting the appropriate $f$ from the previous lemma. 

    Thus the same argument as above shows that each kernel function $k_x$ of the restricted kernel, with $x\in b$, is linearly independent, as the subspace is still Hausdorff. 

    To prove the latter part of the proposition, let us assume we have kernels on $X\times Y\times X\times Z\times W$, $X\times Y$, and $X$ so that $k_{X\times Y\times X\times Z\times W}|_{(X, Y, 0,0, 0)}\cong k_{X\times Y}$, and \[k_{X\times Y\times X\times Z\times W}|_{(X, 0, 0,0, 0)} \cong k_{X\times Y\times X\times Z\times W}|_{(0, 0, X,0, 0)} \cong k_X.\]

    Then we define $k_{X\times Z} := k_{X\times Y\times X\times Z\times W}|_{(0,0,X,Z,0)}.$ Then the restriction of $K_{X\times Z}$ to $X$ will be $k_{X\times Y\times X\times Z\times W_{(0,0, X, 0,0)}}$, which is what we require. 
\end{proof}

\begin{proof}[Proof of Theorem \ref{Tychonoff}]
    
    Consider the set $\mathcal A$ of atomic spaces of $\mathcal T$, and select one such space, call it $A$. Since Tychonoff spaces are regular Hausdorff spaces, we can put a strictly positive definite continuous kernel on $A$.

    Now generate the products of $A$, $\prod^nA$. We give them appropriate kernels, as in lemma \ref{products}.

    Label the set of these spaces $A^\times$. Now if any $B \in A^\times$ satisfies an equation \[B = \prod^\infty\Gamma_j,\] where each $\Gamma_j \in \mathcal A$, then we put kernels on all its constituent spaces by restriction, both the atomic spaces and the product spaces. After we have resolved all such equations, we call the new set of spaces $\overline A$. 
    
    If this accounts for all the atomic spaces, then we are done. Otherwise, we choose a space $A_2$, give it a kernel, and take its products with itself and all the spaces in $\overline A$, repeating the whole process above. In the end we get a set $\overline A_{12}$. After Sysiphean repetition, we achieve our aim. 
\end{proof}

If our category of topological spaces is not decomposable, we can approximate it by choosing termination points of every sequence of projections, and removing everything that would come after. In this way you choose a set of atomic spaces that allows you to put a universal kernel on everything ``above" the chosen atoms.

Indeed, kernels can be put on the spaces that produce the atoms, allowing infinite regression, but whether they can be made to satisfy the axioms of a universal family of kernels is unknown.

A couple of useful corollaries follow from this result. 

\begin{corollary}
    If $\mathcal T$ is a category of second countable compact Tychonoff spaces, then the topology is faithfully mapped by the functor $\mathcal K$. 
\end{corollary}
\begin{proof}
    Each topological space $X \in ob\mathcal T$ is embedded within $\mathcal K(X)$ via $x\mapsto k_x$. Each continuous function $\varphi: X \rightarrow Y$ is faithfully represented via $C^*_\varphi k_x = k_{\varphi(x)}$. These are true without compactness. 

    If $\mathcal T$ is a category of compact spaces, then $K$ applied to each space $T \in ob\mathcal T$ is a continuous embedding of a compact space, and so restriction to the image $\mathcal K(T)$ (considered as a function) provides a homeomorphism. 
\end{proof}
In general, the continuity of the kernel guarantees that the image topology will not be finer than the original topology.

\begin{corollary}
    If $\mathcal T$ is a decomposable category of second countable compact Hausdorff spaces, then the universal family of kernels can additionally be universal kernels.

    That is, the RKHSs provide pointwise approximations of all continuous functions on each space $X \in ob\mathcal T$. 
\end{corollary}
\begin{proof}
    Every compact Hausdorff space is Tychonoff. 

    Let $\mathcal B$ be the basis of the topology $T$. Then as in the proof of lemma \ref{regular kernel}, we construct a kernel function. However, instead of using that as our kernel function, we note that it is a compact subset of $\ell^2$, as it is an injection from a compact space.
    
    We then apply one of several universal kernels on compact subsets of $\ell^2$, such as the radial or Taylor kernels described in [Tan25]. 
    
    The lemmata \ref{products} and \ref{restrictions} work in the same way; the product of universal kernels on spaces $X$ and $Y$ is itself universal on $X\times Y$, and by the definition any universal kernel will be universal on a restriction of its space. 
\end{proof}

The use of a universal kernel on compact $V \subset \ell^2$ is used in [Tan25] to find universal kernels on compact metric spaces. 

\subsection{Reproducing kernels on non-second countable spaces}
Let us consider reproducing kernels on general Tychonoff spaces. 

Following the argument laid out in Lemma \ref{regular kernel}, we can construct a continuous kernel on any Tychonoff space $T$. We do this using functions that are positive at our selected point and zero outside of a selected neighborhood, but we cannot guarantee its behavior on the rest of the set. Thus we construct a continuous kernel, but we cannot guarantee that it is strictly positive definite.

Choose a countable set of points $x_i \in T$. Because $T$ is Tychonoff, we can choose an open neighborhood $U_{x_i}$ of $x_i$, as well as a function $f_i$ so that $f_i(x_i) = 1$ and $f_i|_{T \setminus U_{x_i}} \equiv 0$. 

Indeed, we can choose additional open sets and functions so that for each $x_j \in \{x_i\},$ $x_i \not = x_j$, we have an open set $U_{i\setminus j}$ which contains $x_i$ but not $x_j$. Then for a final step in constructing families of open sets around each $x_i$, we take all finite intersections of $U_{i\setminus j}$ varying over $j$. 

Thus we construct a family $\mathcal U_i$ of open sets around $x_i$, so that $x_i$ is disjoint from any finite set $\{x_j\}_{j\in J}$, $i\notin J$. Each open set $u \in \mathcal U_i$ admits a function $f_u$ so that $f_u(x_i) =1$ and $f_u|_{T\setminus u}\equiv 0$. 

The family $\mathcal U_i$ is countable. Since we have selected only countably many such points $x_i$, the set  $\mathcal U = \bigcup_i\mathcal U_i$ is also countable, and so we can create a continuous kernel by \[k(x, y)= \sum_{u \in \mathcal U} \alpha_{u}f_{u}(x)f_{u}(y)\]
where $\sum_u\alpha_u = 1$. 

The set $k_{x_i}$ is linearly independent, because for any finite set $\{x_i\}_{i\in I}$ there is an open set $u_i\in \mathcal U_i$ so that $x_i\in u_i$ but $x_j\notin u_i$ for all $i\not= j\in I$. 

Then the arguments laid out in lemmata \ref{products} and \ref{restrictions} follow with only aesthetic changes, considering the coverings $\mathcal U$ of the chosen $x_i$ rather than open coverings of the whole space. 

\begin{remark}
    Note that this does not guarantee that a given point $y$ in the neighborhoods of the $x_i$ will be approximated. They may just so happen to be ``missed" by every function $f_u$ where $y \in u$. 
\end{remark}

If we want to guarantee that certain points in product spaces are given nonzero kernel functions, then under the projection to the atomic spaces one can guarantee that the constituent elements are given nonzero kernel functions. The countability requirement does not break in this finite restriction. 

On the other hand, if $\mathcal T$ is locally second countable, then locally we can make all the same arguments as in the proof of Theorem \ref{Tychonoff}, in which case we can guarantee that such points $y$ will not be missed. 

Thus we have proven the following. 

\begin{proposition}
    Let $\mathcal T$ be a decomposable category of Tychonoff spaces. 
    \begin{enumerate}
        \item We can construct a universal family of continuous kernels for $\mathcal T$. 
        \item The family of kernels can be chosen to be strictly positive definite on any countable subset of each $T\in ob\mathcal T$. 
        \item If every $T\in ob\mathcal T$ is locally second countable, the family can be made strictly positive definite on certain open sets.  
    \end{enumerate}
\end{proposition}

If $T$ is both locally second countable and has a countable dense set, then it is second countable and we return to the previous section. Thus if we lack second countability we must expect some tradeoff; either we can faithfully approximate neighborhoods of points, or we can approximate dense subsets of points.

\section{free universal monodromy}
\subsection{Free functions}
Let $Q$ be a quiver, that is a directed multigraph. We consider the \textit{free category}, also known as the \textit{path category}, whose objects are vertices of $Q$ and morphisms are paths (written from right to left) in $Q$. Then if $C$ is an additive category, such as a category of vector spaces with linear transformations between them, we define the functor category $C^Q.$ 

Let $W_{u,v}$ be the set of paths from $u$ to $v$ in $Q$. A free polynomial $\mathbb C\langle Q\rangle(u,v)$ consists of a finite sum \[f= \sum_{w\in W_{u,v}}a_ww.\] If we have $f \in \mathbb C\langle Q\rangle(u,v)$ and $g = \sum_pb_pp \in \mathbb C\langle Q\rangle(v,r)$ we define \[g\cdot f = \sum_w\sum_pb_pa_wp\cdot w \in \mathbb C\langle Q\rangle(u,r)\], where multiplication is path concatenation and thus each $w\cdot p \in W_{u,r}$. 

These polynomials can be evaluated on representations $X \in C^Q$, by replacing each $w$ with the representation $X(w)$. If $C$ is the category of finite dimensional vector spaces, these representations will consist of matrices, and so the evaluations will be polynomials in matrix variables. 

These so-called free polynomials are the archetypal \textit{free functions,} which are defined by a key characteristic of free polynomials; if there is a natural transformation $\Gamma$ satisfying $X \Gamma = \Gamma Y$, then $f(X)\Gamma = \Gamma f(Y).$ 

\begin{definition}
    A subcategory $C$ of $D$ is said to be an \textit{induced subcategory} if whenever $u, v \in ob C$ then $a: u\rightarrow v$ is in $mor C$ if $a\in morD$. That is, all the eligible morphisms are included. 

    An induced subcategory $U$ of the functor category $C^Q$ is said to be a \textbf{full nc-subcategory} of $C^Q$ so long as it is closed under direct sum and natural automorphism. 
    That is: if $X, Y \in U$, then $X\oplus Y \in U$; if $S$ is a natural isomorphism then $S^{-1}XS \in U$. 

    A \textbf{free map} is a function $f$ whose domain is an nc-subcategory $U$ of $C^Q$ and whose codomain is $C^R$, for some free path categories $Q, R$ with the same number of vertices and additive category $C$ if for every $X, Y \in U$ if there is a natural transformation $\Gamma$ in $U$ so that \[X\Gamma = \Gamma Y\] we also have \[f(X)\Gamma = \Gamma f(Y). \]
\end{definition}

We can alternatively characterize a free map as a map that preserves direct sums and natural isomorphisms. 

\begin{definition}
    Let $C^Q$ be a functor category, with $Q$ a free category and $C$ additive. Then if $V: ob Q \rightarrow ob C$ is any function, we define $C^Q_V$ as the the induced subcategory containing only the functors $X \in C^Q$ so that $X(u) = V(u)$ for each $u \in ob Q$. This is not generally an nc-subcategory. 

    In the typical setup, a free map $f$ maps $X \in C^Q_V$ to $f(X) \in C^R_V,$ that is it is graded or ``preserves levels." 
\end{definition}

\subsection{Automatic analyticity}

For our purposes, we shall assume that there is a metric space structure on each homset of our additive category $C$. Then if $Q$ has finitely many objects, we can impose the product metric on each level $C^Q_V$ of $C^Q$.

Then a path is defined as a smooth function $\gamma: (0,1) \rightarrow C^Q_V$, as is typical.

\begin{definition}
    Suppose we have free open set $U \subset C^Q$, and let $f: U\rightarrow C^R$ be a locally bounded free function.

    Then if there is an open set $V \supset U$ and a locally bounded free function $F: V \rightarrow C^R$ with $F|_{U} = f$, then we call $F$ a \textbf{free bornologic continuation} of $f$. 

    We say that $f$ \textbf{freely continues} on $V$.
\end{definition}

This is clearly an analogue of the notion of analytic continuation. Indeed, it is a proper generalization, due to the following foundational result of free analysis [KVV14, AM18].

\begin{proposition}
    Let $C$ be the category of finite dimensional complex vector spaces. 

    Then $f: C^{\langle x_1, \dots, x_d\rangle} \rightarrow C^{\langle x_1, \dots, x_e\rangle}$ is a locally bounded free map if and only if it is analytic. 
\end{proposition}

We can prove a similar theorem in the more general case. An important difference to note is that a robust connection between differentiability and representation as power series has not been established in all additive categories.

\begin{proposition}
    Let $C$ be an additive category whose homsets are metric spaces, and let $Q, R$ be small categories. For a free subcategory $U \subset C^Q$, let $f: U \rightarrow C^R$ be a free map. Then the following are equivalent.
    \begin{itemize}
        \item $f$ is locally bounded;
        \item $f$ is continuous; 
        \item $f$ is differentiable.
    \end{itemize}
\end{proposition}

Differentiability is well known to imply continuity and thus local boundedness, we will prove the that locally bounded implies continuous, and continuous implies differentiable. 

First we state a lemma that appears multiple times in this proof, as well as in spirit in later proofs. 

\begin{lemma}
If the morphisms in the category $\mathcal X \in C^{Q}$ are of the form $x = \m{X & X\Gamma - \Gamma Y\\& Y}$ then \[f(x) = \M{f(X) & f(X)\Gamma - \Gamma f(Y)\\& f(Y)}\]
\end{lemma}
\begin{proof}\label{corner}
    Observe that $x$ is naturally automorphic to $\tilde x$ where $\tilde x = \m{X &\\&Y}$ by way of the natural transformation $s= \m{1 & \Gamma\\& 1}$. Evaluating $f(x) = s^{-1}f(\tilde x)s$ establishes the lemma.
\end{proof}

\begin{proof}[Local boundedness implies continuity]

    Each evaluation level $C^Q_V$ is a metric space under a product metric. 
    
    We will first prove that locally bounded free functions are continuous. To this end, consider a point $X \oplus X \in U_V$; choose a neighborhood  $B(X\oplus X, r)$  of radius $r$ around this point, with $r$ small enough that  $B(X\oplus X, r)\subset U_V$. In particular, choose $r$ small enough that the closure is contained within $U_V$. 
    
    Then since the closure $\overline {B(X\oplus X, r)}$ is compact and $f$ is locally bounded, we have that for some $k>0$ for every $z\in B(X\oplus X, r)$ $\|f(z)\| < k$.

    Now let $\epsilon > 0$ be given, and consider $B(f(X), \epsilon)$. Let $\delta < \min{r/2, r\epsilon/2C}$ so that $B(X, \delta) \subset U_V$. 

    Now choose $Y \in B(X,\delta)$. Then both $\|Y-X\| < \frac{r}{2}$ and $\|\frac{k}{\epsilon}(Y-X)\| < \frac{r}{2}$. Then by the triangle inequality \[\biggl{\|}\M{Y&\frac{k}{\epsilon}(Y-X)\\&X}-\M{X&\\&X}\biggr{\|} < r.\]

    Thus $\M{Y&\frac{k}{\epsilon}(Y-X)\\&X} \in B(X\oplus X, r),$ so \[\biggl\|f\M{Y&\frac{k}{\epsilon}(Y-X)\\&X}\biggr\| < k.\] 

    But since by lemma \ref{corner} \[f\M{Y&\frac{k}{\epsilon}(Y-X)\\&X} = \M{f(Y) & \frac{k}{\epsilon}(f(Y)-f(X)) \\ &f(X)},\] 

    we conclude that \[\|\frac{k}{\epsilon}(f(Y)-f(X))\| < k,\] that is $f(Y) \in B(f(X), \epsilon).$ Thus, by Cauchy's criterion, $f$ is continuous. \end{proof}

\begin{proof}[Continuity implies differentiability]
    Let $f: U \rightarrow C^R$ be continuous. 

    Choose $X \in U$. Then $X \oplus X\in U$. Because $U$ is open in a metric space, we can choose sufficiently small $H$ so that $\m{X & H\\&X} \in U$. Similarly, we can choose a scalar $z$ small enough that \[Z(z) = \M{X +zH & H\\& X} \in U.\]

    If we choose a natural transformation $\Gamma = \frac{1}{z} I$, then we have that in fact \[Z(z) = \M{X+zH & (X +zH)\Gamma - \Gamma X \\& X}.\]
    Our lemma then tells us that \[f(Z(z)) = \M{f(X+zH) & \frac{f(X+zH)-f(X)}{z} \\& f(X)}.\]
    Since $f$ is continuous, we can take the limit and attain the Fr\'echet directional derivative in the top right corner.

\end{proof}

Note that the equivalence of continuity and differentiability does not require a metric space structure; continuity implies convergence of sequences generally, which allows us to calculate the derivative as shown above. We require a metric as ``locally bounded" is not very meaningful outside of metric spaces. 

The following theorem and proof are a natural lifting of the proof in [Pas26], and is reproduced here for completion. 

\begin{theorem}[Universal Free Monodromy]
\label{monodromy}
    Suppose $f$ is continuous on an nc domain $D \subset C^Q_V$. If $f$ freely continues along every path in the nc domain $U \subset C^Q$, then it freely continues to all of $U$. 
\end{theorem}
\begin{proof}
    Suppose $f$ freely continues along the two paths $\gamma_1, \gamma_2$ in $C^Q$ that share the same endpoints, say $s$ and $t$, parametrized along $(0,1)$. Then it must also freely continue along the path \[\gamma_1\oplus \gamma_2 = \M{\gamma_1 &\\&\gamma_2} \in C^Q\] from $s \oplus s$ to $t\oplus t$. 

    For sufficiently small $\varepsilon > 0$, then, the same free bornologic continuation of $f$ extends along the path \[\M{\gamma_1 & \varepsilon\frac{\gamma_1-\gamma_2}{\|\gamma_1 - \gamma_2\|^{1/2}}\\&\gamma_2} = \M{1 & \frac{\varepsilon}{\|\gamma_1 - \gamma_2\|^{1/2}}\\&1}^{-1}
    \M{\gamma_1 &\\&\gamma_2}
    \M{1 & \frac{\varepsilon}{\|\gamma_1 - \gamma_2\|^{1/2}}\\&1}\].

    Because the free bornologic continuation of $f$  is free, we have that

    \[f\M{\gamma_1 & \varepsilon\frac{\gamma_1-\gamma_2}{\|\gamma_1 - \gamma_2\|^{1/2}}\\&\gamma_2} = \M{1 & \frac{\varepsilon}{\|\gamma_1 - \gamma_2\|^{1/2}}\\&1}^{-1}
    f\M{\gamma_1 &\\&\gamma_2}
    \M{1 & \frac{\varepsilon}{\|\gamma_1 - \gamma_2\|^{1/2}}\\&1},\] that is 

    \[f\M{\gamma_1 & \varepsilon\frac{\gamma_1-\gamma_2}{\|\gamma_1 - \gamma_2\|^{1/2}}\\&\gamma_2} = \M{f(\gamma_1) & \varepsilon\frac{f(\gamma_1)-f(\gamma_2)}{\|\gamma_1 - \gamma_2\|^{1/2}}\\&f(\gamma_2)}.\]

    Let us suppose the free bornologic continuation disagreed; at the end points this evaluation would blow up. But this contradicts the local boundedness of the extension along our shifted path, and thus along our original paths as well.
\end{proof}

\section{explaining overfitting}

Now we turn our attention to the stated purpose; providing a theoretical explanation for the robustness of neural networks against over-parametrization. 

Let us now suppose we wish to optimize a neural network $f: \mathbb B_n \rightarrow \mathbb B_m$ with respect to a certain metric. We will model this with a normal state, which is a normalized positive semidefinite function on operators. We do this because states provide a straightforward and robust way of measuring operators, and thus allow us to effectively measure the underlying function $\varphi$. 

\begin{definition}
    A \textbf{state} $\rho$ on an algebra of operators $\mathcal A = \mathcal B(\mathcal H)$ is a linear functional so that $\rho(1)=1$ and $\rho(a^*a) \geq 0$ for every $a \in \mathcal A$.
    
    A \textbf{vector state} $\rho_x$ is given by $\rho_x(a) = \frac{x^*ax}{\|x\|^2}$ for $a\in \mathcal A$, $x\in \mathcal H$. 

    A state $\rho$ is said to be \textbf{normal} if it is weak operator continuous on the unit ball of $\mathcal A$ (or, equivalently, many other things) [KR86, Theorem 7.1.12].

    A normal state is a convex sum of vector states; that is, if $\rho$ is a normal state then \[\rho = \sum_{n \in \mathbb N}\lambda_n\rho_{x_n}\] where $\sum_{n\in\mathbb N}\lambda_n=1$ and $\rho_{x_n}$ is a vector state with vector $x_n \in \mathcal H$. 
\end{definition}
In order to apply our states effectively, we also define the following.

\begin{definition}
    We say a Koopman functor $\mathcal K: \mathcal T \rightarrow RKHS$ is a \textbf{stable model} of $\mathcal T$ if for every continuous map $\varphi$ within $\mathcal T$ the operator $C_\varphi$ within $\mathcal K(\mathcal T)$ is bounded.
    Alternatively, for a class of continuous maps $\{\psi\} \subset mor\mathcal T$ we say that $\mathcal K$ \textbf{stably models $\{\psi\}$} if each $C_\psi$ is a bounded operator within $\mathcal K(\mathcal T)$. 

    Note that if $\mathcal K$ stably models a class $\{\psi\}$, then it stably models the class of such $\psi$ and their composites, giving us a subcategory $\overline{\{\psi\}} $ of $\mathcal T$ which is stably modeled by $\mathcal K$. 
\end{definition}
\begin{example}
    Let $\mathfrak B$ be the category whose objects are the unit balls $\mathbb B_n$ for all $n\in \mathbb N$, and whose morphisms are holomorphic maps $\varphi: \mathbb B_n \rightarrow \mathbb B_m$ such that \[\biggl (\frac{1-\langle \varphi(z), \varphi(w)\rangle}{1-\langle z, w\rangle}\biggr)^n\] is a kernel function on $\mathbb B_n$. 

    The universal family of Drury-Arveson kernels gives us a stable model of $\mathfrak B$. 

    This is a result of theorem 4.2 of [KM25]. 
\end{example}
\begin{example}
    Let $\mathfrak R$ be the category of topological spaces $\mathbb R^n$, with the morphisms containing all analytic maps. 

    Then nearly any kernel we pick will not stably model all morphisms. As shown in [IIS22], for large classes of positive definite kernel functions on $\mathbb R^d$, the only maps that induce bounded composition operators are affine maps. 
\end{example}

Since optimizing over unbounded operators is subtle, in the sequel we will consider only stable models. 

Then we choose a normal state $\rho = \sum_n\lambda_n\rho_{f_n}$, where each $f_n$ is in the RKHS $ \mathcal K(X)$. 

Then if $\varphi \in morX$, the composition operator $C_\varphi$ is bounded and we choose each $f_n$ so that $\|f_n\|=1$, then  
\begin{equation}
     \rho C_\varphi = \sum_n\lambda_nf_n^*C_\varphi f_n.
\end{equation}
\begin{lemma}
    When considered as a function in the variables $f^*, f, C_\varphi$, the map $F(f^*, f, C_\varphi) = \rho C_\varphi$ is a free function (as in [B25]), but can be viewed as a functional. 
\end{lemma}
\begin{proof}
    Allowing direct sums of operators on $\mathcal K(X)$ as well as direct sums of the RKHS elements $f$ gives us \[\begin{split}
    F(f\oplus g, f^* \oplus g^*, C_\varphi \oplus C_\psi) &= \sum_n\lambda_n(f_n^*\oplus g^*_n)(C_\varphi \oplus C_\psi)( f_n\oplus g_n)\\&= \sum_n\lambda_n\biggl(f_n^*C_\varphi f_n\oplus g_n^*C_\psi g_n\biggr).
\end{split}\] Thus it respects the direct sum structure, and the polynomial nature of it respects similarity as well, and is therefore a free function.

If we consider the direct sum spaces as internal direct sums of RKHSs, this becomes \[\begin{split}
    \sum_n \lambda_n (f^*_n+g^*_n)(C_\varphi f_n+C_\psi g_n) &=  \sum_n\lambda_n(f^*_nC_\varphi f_n + g^*_nC_\varphi f_n + f^*_nC_\psi g_n + g^*_nC_\psi g_n)\\ &= \lambda_n\biggl(f_n^*C_\varphi f_n+ g_n^*C_\psi g_n\biggr)
\end{split}\] 
because $f^*C_\psi g$ is calculated as $\M{f&0}^*\m{0 &0\\0&C_\psi}\M{0 & g} =0,$ as in Proposition \ref{between spaces}.
\end{proof} \begin{remark}
    Note that $F$ is not still a state under the direct sum, even if both $\|f_n\|=1=\|g_n\|.$ If, however, we have a convex combination of $f$ and $g$, it is again a state. 

    Thus, positive functionals are free functions. To get a state, we look at a special case. 
\end{remark}

This can be realized as an equation of functions; in particular, we can use the definition of the RKHS to show that \begin{equation}
    \rho C_\varphi = \sum_n\lambda_n\sum_{{i_n}}a_{i_n}f_n({\varphi(x_{i_n})})
\end{equation}
where $f_n = \sum_{i_n}a_{i_n}k_{x_{i_n}}.$ 

Then we wish to optimize $\rho C_\varphi$ in terms of $\varphi$. If we aim to minimize $\rho C_\varphi$ directly we may run into the issue that $\sup_\varphi \|C_\varphi\| =\infty$ despite the operators being bounded, so we will instead aim to minimize $\frac{1}{\|C_\varphi\|}\rho C_\varphi$, which is bounded in norm by $\|\rho\| =1$. 

Then given an arbitrary map $\varphi$, we optimize using standard methods to arrive at a local minimum, namely $\psi$. Since $\rho C_\varphi$ is a free function, Theorem \ref{monodromy} tells us that the extension of $\rho C_\varphi$ to include $\rho C_\psi$ is unique. 

This gives us a local minimum, but in no way should we assume it is a global minimum. In order to find a global optimal, we will essentially compute optimal values in parallel. 

Let us randomly select continuous maps $\varphi_i$ for $i\in [n]=\{0, \dots, n\}$, selected by a distribution $\mu$. Then we choose a state $\rho_{n,\eta}$  given by a convex combination of the $C_{\varphi_i}$ evaluated on the state $\rho$. That is, 
\begin{equation}
    \rho_{n, \eta}\bigoplus_{i\in [n]}C_{\varphi_i} = \sum_i\frac{\eta_i}{\|C_{\varphi_i}\|}\cdot\rho C_{\varphi_i}
\end{equation}
where $\sum_i\eta_i=1$. 

If each of these $\varphi_i$ is individually optimized to get a locally optimal $\psi_i$, then we can find $\inf_\eta\rho_{n,\eta}\bigoplus_{i\in [n]}C_{\psi_i}$. This in effect chooses the best $\psi_i$, and puts all the weight of $\eta$ onto the associated $\eta_i$. 

\begin{remark}
    The set of $\varphi_i$ can be considered a set of parallel processes searching the space. 
\end{remark}

As $n\rightarrow \infty$, so long as the distribution $\mu$ does not miss whole swaths of the space of functions $\phi : X\rightarrow X$, we have that \[ \inf_\eta\rho_{n,\eta}\bigoplus_{i\in [n]}C_{\psi_i} \rightarrow  \inf_\varphi\frac{1}{\|C_\varphi\|}\rho C_\varphi.\]

\begin{remark}[Interpreting beneficial overfitting]
    This can be interpreted in the context of neural networks in the following way. 
    
    In a sparse initialization of a neural network, the sparsity effectively causes several independent neural nets to be set up in parallel. Then, over the course of optimization, each of these neural nets is individually optimized, as is the weighting of the neural nets in parallel. 

    As such, when the number of nodes is increased, we are in effect adding additional neural nets which will ultimately be optimized and selected between. 

    Thus increasing the number of nodes opens up the possibility space searched to find a solution.
\end{remark}

This same framework can be leveraged more broadly, in the case where we are combining multiple simple ``utility" functions into a broader tool.

\begin{remark}
    If we have a finite collection $\Psi$ of optimized ``utility" neural networks mapping $X$ to itself, we choose our tuple $x_i = C_{\psi_i}$ for $ \psi_i \in \Psi.$ This tuple $x$ is a point in $\mathcal B^n$, where $\mathcal B$ the vector space of composition operators generated by $\Phi = \{\varphi : X\rightarrow X\}$. Thus we have a $\varepsilon$-neighborhood of $x$ within $\mathcal B$, because we are assuming that our functions are stably modeled. 

    Now suppose we have an objective functional $\rho$, that we wish to optimize against using our tools. Since we are looking for the best combination of these functions, we will cap the weights of these functions. For simplicity, let us bound the weights so that the coefficients are non-negative and must sum to $1$. Thus our functions are of the form \[f = \sum_wa_ww,\] where $w$ is a multi-index. Since this is a bounded convex space of functions, we are guaranteed a (not necessarily unique) optimum, say $\tilde f$, so that $\rho\tilde f \geq \rho f.$ Since this function is analytic, it can be uniquely extended to the neighborhood of $x$. 
\end{remark}

\newpage

\section{References}

[AM18] Agler, Jim and John E. McCarthy. ``Global holomorphic functions in several noncommuting variables." \textit{Canadian journal of mathematics} 67, issue 2, pp 241-285. (2018)

[B25] Bushelli, Julian. ``The free functional calculus in general." arXiv:2506.00170 (2025)

[EEM22] Elefante, Giacomo; Erb, Wolfgang; Marchetti, Francesco; Perracchione, Emma; 
Poggiali, Davide; and Gabriele Santin. ``Interpolation with polynomial kernels." 	arXiv:2212.07658 (2022)

[HKM12] Helton, J. William; Klep, Igor; McCullough, Scott. ``The convex positivstellensatz in a free algebra." \textit{Advances in mathematics,}
 Volume 231, 1, pp. 516-534 (2012)

[IIS22] Ikeda, Masahiro; Ishakawa, Isao; and Yoshihiro Sawano. ``Boundedness of composition operators on reproducing kernel Hilbert spaces with analytic positive definite functions." \textit{Journal of mathematical analysis and applications,} 511, issue 1. (2022)

[JMS21] Jury, Michael; Martin, Robert T. W.; Shamovich, Eli. ``Non-commutative rational functions in the full fock space." \textit{Transactions of the American mathematical society}, Vol. 374, 9 (2021)

[KM25] Kumari, Preeti and P. Muthukumar. ``Composition operators on reproducing kernel Hilbert spaces,"  arXiv:2509.14655v1 (2025)

[Koo31] Koopman, B. O. ``Hamiltonian systems and transformations in Hilbert space." \textit{Proceedings of the national academy of sciences} 17 issue 5, pp 315-318. (1931)

[KR86] \textit{Fundamentals of the theory of operator algebras II,} Richard V. Kadison and John R. Ringrose. AMS. (1986)

[KVV14] Kalyuzhnyi-Verbovetski\i, Dmitry S., and Victor Vinnikov. ``Foundations of noncommutative function theory." Mathematical surveys and monographs; volume 199. (AMS, 2014) 

[Mun00] \textit{Topology,} James R. Munkres. Pearson, 2000. 

[Pas23] Pascoe, James Eldred. ``Invariant structure preserving functions and an Oka-Weil Kaplansky density type theorem." Tohoku (2023): 465-482.

[Pas26] \textit{Noncommutative plurisubharmonicity, existence of plurisubharmonic conjugates, and free universal monodromy,} James E. Pascoe. Transactions of the American Mathematical Society, 20 Feb. 2026. 

[Rie17] Riehl, Emily. \textit{Category theory in context}. Courier Dover Publications, 2017.

[SaSa16] \textit{Theory of Reproducing Kernels and Applications}, Saburou Saitoh and Yoshihiro Sawano. Springer. (2016)

[Tan25] \textit{Explicit Universal and Approximate-Universal Kernels on Compact Metric Spaces}, Eloi Tanguy, 2025. HAL Id: hal-05084408.

[Tay72] Taylor, Joseph L. ``A general framework for a multi-operator functional calculus." \textit{Advances in mathematics} 9, pp. 183-252 (1972)

[ZBHRO21] Zhang, Chiyuan; Bengio, Samy; Hardt, Moritz; Recht, Benjamin; and Oriol Vinyols. ``Understanding deep learning (still) requires rethinking generalization." \textit{Communications of the ACM,} 64 issue 3, pp 107-115. (2021)

[Zwi08] Zwicknagel, Barbara. ``Power Series Kernels", \textit{Journal of Constructive Approximation.} 29, pp 61-84. (2008)

\end{document}